\documentclass[12pt,a4paper,twoside]{article}
\usepackage[french] {babel}
\usepackage{amssymb,amsmath,amsfonts,graphics}
\usepackage{stmaryrd}
\usepackage[all]{xy}

\begin{document}

\newcommand{\hooklongrightarrow}{\lhook\joinrel\longrightarrow}

\begin{center}
\textbf{Cohomologie syntomique:\\
 liens avec les cohomologies \'etale et rigide}

\vskip30mm

Jean-Yves ETESSE
  \footnote{(CNRS - Institut de Math\'ematique, Universit\'e de Rennes 1, Campus de Beaulieu - 35042 RENNES Cedex France)\\
E-mail : Jean-Yves.Etesse@univ-rennes1.fr}
\end{center}

\vskip50mm
\noindent\textbf{Sommaire}\\

		\begin{enumerate}
		\item[] Introduction
		\item[1.] Site syntomique
		\item[2.] Cohomologie syntomique \`a supports compacts
		\item[3.] Comparaison avec cohomologie \'etale et  cohomologie rigide
	\end{enumerate}

\newpage

\noindent\textbf{R\'esum\'e}\\

La cohomologie syntomique ici d\'efinie fait le lien entre la cohomologie rigide et la cohomologie \'etale, en interpr\'etant cette derni\`ere comme points fixes du Frobenius agissant sur la premi\`ere.

Soit $\mathcal{V}$ un anneau de valuation discr\`ete complet, de corps r\'esiduel parfait $k = \mathcal{V}/\mathfrak{m}$ de caract\'eristique $p > 0$ et de corps des fractions $K$ de caract\'eristique $0$. Apr\`es avoir d\'efini la cohomologie syntomique \`a supports compacts d'un faisceau ab\'elien $\mathcal{F}$ sur un $k$-sch\'ema $X$, nous montrons que celle-ci co\"{\i}ncide avec la cohomologie \'etale \`a supports compacts lorsque $\mathcal{F}$ est un faisceau lisse. Si de plus le $F$-isocristal convergent associ\'e au faisceau lisse $\mathcal{F}$  provient d'un isocristal surconvergent $E$, alors la cohomologie rigide de $E$ s'exprime elle aussi comme limite de cohomologies syntomiques: la cohomologie \'etale \`a supports compacts de $\mathcal{F}$ est alors les points fixes du Frobenius agissant sur la cohomologie rigide de $E$.\\

\noindent\textbf{Abstract}\\

Syntomic cohomology here defined yields a link between rigid cohomology and etale cohomology, viewing the last one as the fixed points under Frobenius of the former one.

Let $\mathcal{V}$ be a complete discrete valuation ring, with perfect residue field $k = \mathcal{V}/\mathfrak{m}$ of characteristic $p > 0$ and fraction field $K$ of characteristic $0$. Having defined syntomic cohomology with compact supports of an abelian sheaf $\mathcal{F}$ on a $k$-scheme $X$, we show that it coincides with etale cohomology with compact supports when $\mathcal{F}$ is a lisse sheaf. If moreover the convergent $F$-isocrystal associated to $\mathcal{F}$ comes from an overconvergent isocrystal $E$, then the rigid cohomology of $E$ expresses as a limit of syntomic cohomologies: then the etale cohomology with compact supports of $\mathcal{F}$ is the fixed points of Frobenius acting on the rigid cohomology of $E$.\\


Mots cl\'es: cohomologie syntomique, $F$- isocristaux convergents, $F$- isocristaux surconvergents, cohomologie cristalline, cohomologie rigide.\\


\newpage

\section*{ Introduction}

La cohomologie syntomique introduite dans cet article fait le lien entre la cohomologie \'etale et la cohomologie rigide, lien qui sera utilis\'e ult\'erieurement pour r\'esoudre une conjecture de Katz sur les z\'eros et p\^oles unit\'es $p$-adiques des fonctions $L$.\\

Comme pour le topos \'etale [Mi, II, theo 3.10] on obtient au \S1 une description du topos syntomique [th\'eo (1.3)] calqu\'ee sur celle de [SGA 4,T 1, IV, th\'eo (9.5.4)], SGA 4 qui travaille en termes de sous-topos ouvert et du sous-topos ferm\'e compl\'ementaire. On en d\'eduit les suites exactes courtes usuelles de localisation [th\'eo (1.5)].\\

Au \S2 la cohomologie syntomique \`a supports compacts est d\'efinie: en particulier il faut s'assurer de l'ind\'ependance par rapport \`a la compactification choisie [prop (2.1)]. Les suites exactes courtes de localisation du \S1 fournissent alors la suite exacte longue de localisation en cohomologie syntomique \`a supports compacts [th\'eo (2.6)].\\

Au \S3 on \'etablit que les cohomologies syntomique et \'etale \`a supports compacts d'un faisceau lisse $\mathcal{F}$ co\" incident [th\'eo (3.2)]. De m\^eme la cohomologie rigide \`a supports compacts d'un $F$-isocristal surconvergent unit\'e $E$ associ\'e \`a $\mathcal{F}$ co\" incide avec une limite de la  cohomologie syntomique \`a supports compacts d'un $E_{n}^{m\mbox{-}cris}$ associ\'e \`a $\mathcal{F}$ [ th\'eo (3.3.13)(1) et (3.3.16)]. Comme il existe une suite exacte courte sur le site syntomique qui relie  $\mathcal{F}$ et $E_{n}^{m\mbox{-}cris}$ [th\'eo (3.3.12)], on en d\'eduit que la cohomologie \'etale \`a supports compacts de $\mathcal{F}$ s'identifie aux points fixes du Frobenius agissant sur la cohomologie rigide \`a supports compacts de $E$ [th\'eo(3.3.13) (2) et (3)].\\

\section*{1. Site syntomique}

\textbf{1.1.} Si $X$ est un sch\'ema quelconque, le gros site syntomique de $X$ est d\'efini comme suit [SGA 3, IV, 6.3] : la cat\'egorie sous-jacente est celle des sch\'emas sur $X$ et la topologie est engendr\'ee par les familles finies surjectives de morphismes syntomiques (i.e. plats et localement intersection compl\`ete). On rappelle que les morphismes syntomiques sont ouverts, demeurent syntomiques par changement de base, et sont localement relevables le long d'une immersion ferm\'ee. Le gros site (resp. petit site) syntomique de $X$ sera not\'e $SYNT(X)$ (resp. $synt(X)$) et le topos correspondant $X_{SYNT}$ (resp. $X_{synt}$) : lorsqu'on ne voudra pas distinguer entre les deux situations on notera $\mathcal{T}(X)$ (resp. $X_{\mathcal{T}}$) l'un ou l'autre de ces deux sites (resp. topos).\\

Soit $j : U \hookrightarrow X$ une immersion ouverte ; $j$ d\'efinit un couple de foncteurs adjoints $(j_{\ast}, j^{-1})$

$$
U_{\mathcal{T}} \displaystyle \mathop{\rightarrow}^{j^{-1} \atop{\longleftarrow}}_{j _{\ast}} X_{\mathcal{T}}.
$$

\noindent Dans la suite $X_{\mathcal{T}}$ sera annel\'e par un faisceau d'anneaux $\mathcal{A}$ et $U_{\mathcal{T}}$ sera annel\'e par $j^{-1} \mathcal{A}$, not\'e $\mathcal{A}_{\vert U}$ : on note $_{\mathcal{A}}X_{\mathcal{T}}$ (resp. $_{\mathcal{A}\vert U}U_{\mathcal{T}}$) la cat\'egorie des faisceaux des $\mathcal{A}$-modules \`a gauche sur $X_{\mathcal{T}}$ (resp. des $\mathcal{A}_{\vert U}$-modules \`a gauche sur $U_{\mathcal{T}}$). Le foncteur

$$
j^{\ast} : \ _{\mathcal{A}}X_{\mathcal{T}} \longrightarrow\  _{\mathcal{A}\vert U}U_{\mathcal{T}}
$$

\noindent admet un adjoint \`a gauche $j_{!}$ [Mi, II, Rk 3.18] et [SGA 4, IV, $\S$ 14] d\'efini par\\

\noindent
(1.1.1)
$\left\lbrace
\begin{array}{l}
j_{!} (\mathcal{F})(X') = \mathcal{F}(X')\  \textrm{si}\ X' \longrightarrow X\  \textrm{se factorise par}\  U\\
\textrm{et}\  j_{!} (\mathcal{F}) (X') = 0\ \textrm{sinon\ ;}
\end{array}
\right.
$\\

\noindent $j_{!}$ est exact [loc. cit].\\

\noindent Remarquons que $j^{\ast}$ est exact puisqu'il admet un adjoint \`a droite et un adjoint \`a gauche.\\

De m\^eme si $i : Z \hookrightarrow X$ est une immersion ferm\'ee, $i$ d\'efinit un couple de foncteurs adjoints $(i_{\ast}, i^{-1})$ :

$$
Z_{\mathcal{T}} \displaystyle \mathop{\rightarrow}^{i^{-1} \atop{\longleftarrow}}_{i _{\ast}} X_{\mathcal{T}}\ ;
$$

\noindent $Z_{\mathcal{T}}$ sera annel\'e par $i^{-1} \mathcal{A}$, not\'e $\mathcal{A}_{\vert Z}$.\\

\noindent Le foncteur
$$
i^{\ast} :\  _{\mathcal{A}}X_{\textrm{SYNT}} \longrightarrow\ _{\mathcal{A} \vert Z}Z_{\textrm{SYNT}}
$$

\noindent admet un adjoint \`a gauche $i_{!}$ [Mi, II, Rk 3.18] et $i_{!}$ est exact [loc. cit.] : en particulier 
$$i^{\ast} :\  _{\mathcal{A}}X_{\textrm{SYNT}} \longrightarrow\ _{\mathcal{A} \vert Z}Z_{\textrm{SYNT}}$$ est exact.\\
Le foncteur $$i^{\ast} :\  _{\mathcal{A}}X_{\textrm{synt}} \longrightarrow\ _{\mathcal{A} \vert Z}Z_{\textrm{synt}}$$ est lui aussi exact gr\^ace \`a [Mi, II, 2.6 et 3.0 p 68] et [EGA $0_{I}$, 1.4.12] car les morphismes syntomiques demeurent syntomiques par changement de base.\\
De plus le foncteur 
$$i_{\ast} :\  _{\mathcal{A} \vert Z}Z_{\mathcal{T}} \longrightarrow\ _{\mathcal{A}}X_{\mathcal{T}}$$
  est exact, car tout morphisme syntomique se rel\`eve, localement le long d'une immersion ferm\'ee, en un morphisme syntomique.\\

\textbf{1.2.} Soient $i : Z \hookrightarrow X$ une immersion ferm\'ee, et $j : U \hookrightarrow X$ l'immersion ouverte du compl\'ementaire de $Z$.
Pour un $X$-sch\'ema $X'$ on note $U' = X' \times_{X} U$, $Z' = X' \times_{X} Z$; tout recouvrement syntomique $W \twoheadrightarrow Z'$ se rel\`eve localement en $\tilde{W} \rightarrow X'$ syntomique : pour all\'eger l'\'ecriture on supposera le rel\`evement global. Par suite, si $\tilde{U} \twoheadrightarrow U'$ est surjectif syntomique et $\tilde{W}$ tel que ci-dessus, alors $(\tilde{U}, \tilde{W})$ est un recouvrement syntomique de $X'$.

\vskip 3mm
\noindent \textbf{Lemme (1.2.1)}. 
\textit{Avec les notations de (1.2) et pour $\mathcal{F} \in {_{\mathcal{A}}X_{\mathcal{T}}}$ le carr\'e suivant, o\`u les fl\`eches sont les fl\`eches canoniques}

$$
\xymatrix{
\mathcal{F} \ar[r] \ar[d] & j_{\ast} j^{\ast}(\mathcal{F)} \ar[d] \\
i_{\ast} i^{\ast}(\mathcal{F}) \ar[r] & i_{\ast} j_{\ast} j^{\ast}(\mathcal{F})
}
$$
\noindent \textit{est cart\'esien.}

\vskip 3mm
\noindent \textit{D\'emonstration}. On notera $\mathcal{G}$ le produit fibr\'e.\\
Pour un $X$-sch\'ema $X'$ on consid\`ere un recouvrement $(\tilde{U}, \tilde{W})$ de $X'$ du type pr\'ec\'edent : comme $\tilde{U} \rightarrow X'$ et $\tilde{W} \rightarrow X'$ sont deux morphismes syntomiques, on est ramen\'e \`a \'etablir l'isomorphisme $\mathcal{F}\  \tilde{\rightarrow}\  \mathcal{G}$ au-dessus d'un $X'$-sch\'ema syntomique $\tilde{W}$ ; la d\'emonstration se fait donc sur le petit site de $X'$. On pose $W = \tilde{W} \times_{X'}Z'$.\\

Le faisceau $j_{\ast} j^{\ast}(\mathcal{F})$ est le faisceau associ\'e au pr\'efaisceau

$$
\tilde{W} \longmapsto \mathcal{F}(V) \qquad , \quad \textrm{avec}\  V := \tilde{W} \times_{X'} U' ,
$$

\noindent et $i_{\ast} i^{\ast}(\mathcal{F})$ est le faisceau associ\'e au pr\'efaisceau 

$$
\tilde{W} \longmapsto \displaystyle \mathop{\lim}_{\rightarrow \atop{W'}} \mathcal{F}(\mathcal{W}')\ ,
$$
\noindent la limite \'etant prise sur les diagrammes commutatifs

$$\begin{array}{c}
\xymatrix{
\tilde{W} \ar[d] & W' \ar[l] & W \ar[l] \ar[d]\\
X' & &Z' \ar[ll]
}
\end{array}
\leqno{(1.2.2)}
$$

\noindent avec $W' \rightarrow \tilde{W}$ syntomique.\\
De m\^eme $i_{\ast} i^{\ast} j_{\ast} j^{\ast}(\mathcal{F})$ est le faisceau associ\'e au pr\'efaisceau

$$
\tilde{W} \longmapsto \displaystyle \mathop{\lim}_{\rightarrow \atop{W'}} \mathcal{F}(\mathcal{W}' \times_{\tilde{W}} V) =  \displaystyle \mathop{\lim}_{\rightarrow \atop{W'}} \mathcal{F} (\mathcal{W}' \times_{X'} U'),
 $$
 
 \noindent avec $W'$ comme en (1.2.2). On remarque alors que $(V, W')$ est un recouvrement syntomique de $\tilde{W}$ : en effet le $\tilde{W}$-morphisme $W \rightarrow W'$ fournit une section du morphisme syntomique $W' \times_{\tilde{W}} W \rightarrow W$ ; ce dernier est donc surjectif et on conclut comme pour le recouvrement $(\tilde{U}, \tilde{W})$ de $X'$.\\
 Ainsi on a bien un isomorphisme $\mathcal{F}\  \tilde{\rightarrow}\  \mathcal{G}$ au-dessus de $\tilde{\mathcal{W}}$, d'o\`u le lemme. $\square$\\
 
Notons $\mathbb{T}(_{\mathcal{A}}X_{\mathcal{T}})$ la cat\'egorie des triplets $(\mathcal{F}_{1}, \mathcal{F}_{2}, \alpha)$ o\`u $\mathcal{F}_{1} \in _{\mathcal{A} \vert Z}Z_{\mathcal{T}}$, $\mathcal{F}_{2} \in\  _{\mathcal{A} \vert U}U_{\mathcal{T}}$ et $\alpha$ est un morphisme $\alpha : \mathcal{F}_{1} \rightarrow i^{\ast} j_{\ast}  \mathcal{F}_{2}$ ; les morphismes entre deux tels triplets sont d\'efinis de la mani\`ere naturelle, analogue \`a [Mi, II, \S\ 3].
 
 \vskip 3mm
\noindent \textbf{Th\'eor\`eme(1.3)}. 
\textit{Soient $i : Z \hookrightarrow X$ une immersion ferm\'ee de sch\'emas et $j : U \hookrightarrow X$ l'immersion ouverte du compl\'ementaire de $Z$. Le foncteur
$$\mathcal{F} \longmapsto (i^{\ast} \mathcal{F}, j^{\ast} \mathcal{F}, \alpha)$$ o\`u $\alpha$
est le morphisme canonique $\alpha : i^{\ast} \mathcal{F} \rightarrow i^{\ast} j_{\ast} j^{\ast} \mathcal{F}$, induit une \'equivalence de cat\'egories entre $_\mathcal{A}X_{\mathcal{T}}$ et $\mathbb{T}(_\mathcal{A}X_{\mathcal{T}})$.
 }

\vskip 3mm
\noindent \textit{D\'emonstration}. La d\'emonstration est analogue \`a celle de Fontaine-Messing [F-M, 4.4]. Le th\'eor\`eme r\'esulte du lemme (1.2.1) par la m\^eme m\'ethode que pour le site \'etale [Mi, II, theo 3.10]. $\square$\\

En identifiant $_\mathcal{A}X_{\mathcal{T}}$ et $\mathbb{T}(_\mathcal{A}X_{\mathcal{T}})$ via le th\'eor\`eme (1.3) on d\'efinit six foncteurs

$$
\begin{array}{c}
\xymatrix{
& \ar[l]_{i^{\ast}} \quad & \quad  \ar[l]_{j_{!}} &\\
_{\mathcal{A} \vert Z}Z_{\mathcal{T}} \ar[r]^{i_{\ast}} & _{\mathcal{A}}X_{\mathcal{T}} \ar[r]^{j^{\ast}} & _{\mathcal{A}}U_{\mathcal{T}}\\
& \ar[l]_{i^{!}} \quad & \quad  \ar[l]_{j_{\ast}} &
}
\end{array}
\leqno{(1.4)}
$$

\noindent dont la description est la suivante :\\
$i^{\ast} : \mathcal{F}_{1} \leftarrow (\mathcal{F}_{1}, \mathcal{F}_{2}, \alpha : \mathcal{F}_{1} \rightarrow i^{\ast} j_{\ast} \mathcal{F}_{2}),\  j_{!} : (0, \mathcal{F}_{2}, 0) \leftarrow \mathcal{F}_{2}$\\
$i_{\ast} : \mathcal{F}_{1} \mapsto (\mathcal{F}_{1}, 0, 0) \quad \qquad \qquad \qquad ,\ j^{\ast} : (\mathcal{F}_{1}, \mathcal{F}_{2}, \alpha : \mathcal{F}_{1} \rightarrow i^{\ast} j_{\ast} \mathcal{F}_{2}) \mapsto \mathcal{F}_{2}$\\
$i^! : \mbox{Ker} \alpha \leftarrow (\mathcal{F}_{1}, \mathcal{F}_{2}, \alpha : \mathcal{F}_{1} \rightarrow i^{\ast} j_{\ast} \mathcal{F}_{2})$,\\
 $j_{\ast} : (i^{\ast} j_{\ast} \mathcal{F}_{2}, \mathcal{F}_{2}, id : i^{\ast} j_{\ast} \mathcal{F}_{2} \rightarrow i^{\ast} j_{\ast} \mathcal{F}_{2}) \leftarrow \mathcal{F}_{2}. $

\vskip 3mm
\noindent \textbf{Th\'eor\`eme(1.5)}. 
\textit{Avec les notations pr\'ec\'edentes on a :
\begin{enumerate}
\item[(1)] Chaque foncteur est adjoint \`a gauche de celui \'ecrit la ligne au-dessous ; en particulier on a un isomorphisme de transitivit\'e $(j_{1} j_{2})_{!} = j_{1!} j_{2!}$.
\item[(2)] Les foncteurs $i^{\ast}, i_{\ast}, j^{\ast}, j_{!}$ sont exacts ; les foncteurs $j_{\ast}$, $i^{!}$ sont exacts \`a gauche.
\item[(3)] Les compos\'es $i^{\ast} j_{!}, i^{!} j_{!}, i^{!} j_{\ast}, j^{\ast} i_{\ast}$ sont nuls.
\item[(4)] Les foncteurs $i_{\ast}, j_{\ast}$ et $j_{!}$ sont pleinement fid\`eles.
\item[(5)] Les foncteurs $j_{\ast}, j^{\ast}, i^{!}, i_{\ast}$ envoient les injectifs sur les injectifs.
\item[(6)] Pour tout $\mathcal{F} \in _{\mathcal{A}}X_{\mathcal{T}}$  (resp. $\mathcal{F} \in _{\mathcal{A} \vert U}U_{\mathcal{T}})$ on a des suites exactes courtes
\begin{itemize}
\item[(6.1)] $0 \longrightarrow j_{!} j^{\ast} \mathcal{F} \longrightarrow \mathcal{F} \longrightarrow i_{\ast} i^{\ast} \mathcal{F} \longrightarrow 0$
\item[(6.2)] $0 \longrightarrow i_{\ast} i^{!} \mathcal{F} \longrightarrow \mathcal{F} \longrightarrow j_{\ast}j^{\ast} \mathcal{F}$
\item[[resp. (6.3)] $0 \longrightarrow j_{!} \mathcal{F} \longrightarrow j_{\ast} \mathcal{F} \longrightarrow i_{\ast} i^{\ast} j_{\ast} \mathcal{F} \longrightarrow 0].$\\
De plus le couple de foncteurs $(i^{\ast}, j^{\ast}) $ est conservatif.
\end{itemize}
\end{enumerate}
}

\vskip 3mm
\noindent \textit{D\'emonstration}.\par

Le \textit{(1)} et le \textit{(2)} ont d\'ej\`a \'et\'e vus, et l'isomorphisme de transitivit\'e $(j_{1} j_{2})_{!} = j_{1!} j_{2!}$ r\'esulte de la formule $(j_{1} j_{2})^{\ast} = j_{2}^{\ast} j_{1}^{\ast}$.\par

Le \textit{(3)} et le \textit{(4)} r\'esultent des descriptions (1.4).\par

Le \textit{(5)} r\'esulte de \textit{(1)} et \textit{(2)} et du fait qu'un foncteur avec un adjoint \`a gauche exact pr\'eserve les injectifs [Mi, III, 1.2].\par

La suite \textit{(6.3))} provient de \textit{(6.1)} en remarquant que $j^{\ast} j_{\ast} \mathcal{F} = \mathcal{F}.$\par

Compte tenu des identifications (1.4) la suite \textit{(6.1)} s'\'ecrit

$$
0 \longrightarrow (0, j^{\ast} \mathcal{F}, 0) \longrightarrow (i^{\ast} \mathcal{F}, j^{\ast} \mathcal{F}, \alpha) \longrightarrow (i^{\ast} \mathcal{F}, 0, 0) \longrightarrow 0.
$$

\noindent Pour d\'emontrer qu'elle est exacte il nous suffit donc de montrer qu'une suite de $_{\mathcal{A}}X_{\mathcal{T}}$\\

\noindent \textit{(6.4)} $\qquad \qquad  0 \longrightarrow \mathcal{F}' \longrightarrow \mathcal{F} \longrightarrow \mathcal{F}'' \longrightarrow 0$\\

\noindent est exacte si et seulement si les suites\\

\noindent \textit{(6.5)} $\qquad \qquad  0 \longrightarrow i^{\ast}(\mathcal{F}') \longrightarrow i^{\ast}(\mathcal{F}) \longrightarrow i^{\ast}(\mathcal{F}'') \longrightarrow 0$\\

\noindent \textit{(6.6)} $\qquad \qquad  0 \longrightarrow j^{\ast}(\mathcal{F}') \longrightarrow j^{\ast}(\mathcal{F}) \longrightarrow j^{\ast}(\mathcal{F}'') \longrightarrow 0$\\

\noindent de $_{\mathcal{A} \vert Z}Z_{\mathcal{T}}$ et $_{\mathcal{A} \vert U}U_{\mathcal{T}}$ respectivement, sont exactes, i.e. que le couple de foncteurs $(i^{\ast}, j^{\ast})$ est conservatif.\\

L'exactitude de \textit{(6.4)} entra\^{\i}ne celle de \textit{(6.5)} et \textit{(6.6)} puisque $i^{\ast}$ et $j^{\ast}$ sont exacts.\\

R\'eciproquement, l'exactitude de \textit{(6.5)} et \textit{(6.6)} entra\^{\i}ne celle des suites\\

\noindent \textit{(6.7)} $\qquad \qquad  0 \longrightarrow i_{\ast} i^{\ast}(\mathcal{F}') \longrightarrow i_{\ast} i^{\ast}(\mathcal{F}) \displaystyle \mathop{\longrightarrow} ^{\varphi_{Z}} i_{\ast} i^{\ast}(\mathcal{F}'') \longrightarrow 0,$\\

\noindent \textit{(6.8)} $\qquad \qquad  0 \longrightarrow j_{\ast} j^{\ast}(\mathcal{F}') \longrightarrow j_{\ast} j^{\ast}(\mathcal{F}) \displaystyle \mathop{\longrightarrow} ^{\varphi_{U}} j_{\ast} j^{\ast}(\mathcal{F}'') ,$\\

\noindent \textit{(6.9)} $\qquad \qquad  0 \longrightarrow i_{\ast} i^{\ast} j_{\ast} j^{\ast}(\mathcal{F}') \longrightarrow i_{\ast} i^{\ast} j_{\ast} j^{\ast}(\mathcal{F}) \longrightarrow i_{\ast} i^{\ast} j_{\ast} j^{\ast}(\mathcal{F}''),$\\

\noindent d'o\`u l'exactitude de\\

\noindent \textit{(6.10)} $\qquad \qquad \qquad 0 \longrightarrow \mathcal{F}' \longrightarrow \mathcal{F} \longrightarrow \mathcal{F}''$\\

\noindent gr\^ace au lemme (1.2.1).\\

Donc pour tout $\mathcal{F} \in _{\mathcal{A}}X_{\mathcal{T}}$ on a l'exactitude de la suite

$$
0 \longrightarrow j_{!}j^{\ast} \mathcal{F} \longrightarrow \mathcal{F} \longrightarrow i_{\ast} i^{\ast} \mathcal{F}.
$$

\noindent Montrons la surjectivit\'e de $\rho : \mathcal{F} \rightarrow i_{\ast}i^{\ast} \mathcal{F}$. Comme pour le lemme (1.2.1), dont on utilise les notations, on est ramen\'e au petit site de $X'$. Soient $\tilde{W}$ un $X'$-sch\'ema syntomique, $V := \tilde{W} \times_{X'}U'$ et $s \in i_{\ast} i^{\ast}(\mathcal{F})(\tilde{W}) = \displaystyle \mathop{\lim}_{\rightarrow \atop{W'}} \mathcal{F}(W')$, o\`u $W'$ est comme dans (1.2.2) : il existe un $W'$ et $s_{W'} \in \mathcal{F}(W')$ d'image $s$. On a vu dans la preuve du lemme (1.2.1) que $(V, W')$ est un recouvrement syntomique de $\tilde{W}$ : notons $s'$ l'image de $s$ par l'application restriction

$$
i_{\ast} i^{\ast} \mathcal{F}(\tilde{W}) \longrightarrow i_{\ast} i^{\ast} \mathcal{F}(W')
$$

$$
s \mapsto s'\ ;
$$

\noindent alors $s_{W'} \in \mathcal{F}(W')$ a pour image $s'$ par $\rho$.\\

\noindent Comme $i_{\ast} i^{\ast}(\mathcal{F})(V) = 0$, tout $s_{V} \in \mathcal{F}(V)$ est un rel\`evement de $s$ dans $i_{\ast} i^{\ast}(\mathcal{F})(V) = 0$. D'o\`u la surjectivit\'e de $\rho$, et l'exactitude de \textit{(6.1)}.\\

On a alors un diagramme commutatif \`a lignes exactes

$$
\xymatrix{
0 \ar[r]  & j_{!} j^{\ast} \mathcal{F} \ar[r] \ar[d]_{u} & \mathcal{F} \ar[r] \ar[d]_{v} & i_{\ast} i^{\ast} \mathcal{F} \ar[r] \ar[d]_{w} & 0 \\
0 \ar[r] & j_{!} j^{\ast} \mathcal{F}'' \ar[r] & \mathcal{F}'' \ar[r] & i_{\ast} i^{\ast} \mathcal{F}'' \ar[r] & 0.
}
$$

\noindent L'exactitude de \textit{(6.5)} et \textit{(6.6)} et celle des foncteurs $j_{!}$ et $i_{\ast}$ prouve que $u$ et $w$ sont surjectifs : d'o\`u la surjectivit\'e de $v$ ; jointe \`a l'exactitude de \textit{(6.10)} ceci prouve que le couple $(i^{\ast}, j^{\ast})$ est conservatif.\\

\noindent L'exactitude de \textit{(6.2)} en r\'esulte via la description de $i^{!}$ fournie en (1.4). $\square$

\section*{2. Cohomologie syntomique \`a supports compacts}

Soit $X$ un sch\'ema s\'epar\'e de type fini sur un corps $k$ ; on sait par Nagata que $X$ est ouvert dans un $k$-sch\'ema propre $\overline{X}$ ; on note $j : X \hookrightarrow \overline{X}$ l'immersion ouverte.\\

On se place sous les notations de (1.1) : $\overline{X}_{\mathcal{T}}$ est annel\'e par un faisceau d'anneaux $\mathcal{A}$ et $\mathcal{F}$ est un faisceau de $\mathcal{A}$-modules.\\

Si $\mathcal{F}$ est \'element de $_{\mathcal{A}}\overline{X}_{\textrm{synt}}$, on notera encore $\mathcal{F}$ son image inverse dans $_{\mathcal{A}}\overline{X}_{\textrm{SYNT}}$ par le morphisme de topos $\overline{X}_{\textrm{SYNT}} \rightarrow \overline{X}_{\textrm{synt}}$ et pour tout entier $i \geqslant 0$, on a [E-LS 2, (1.2)]

$$
H^i(\overline{X}_{\textrm{synt}}, \mathcal{F}) = H^i(\overline{X}_{\textrm{SYNT}}, \mathcal{F}),
$$
\noindent et de m\^eme pour la topologie \'etale.

\vskip 3mm
\noindent \textbf{Proposition - D\'efinition (2.1)}. \textit{Sous les hypoth\`eses pr\'ec\'edentes, si $\mathcal{A}$ est de torsion et $\mathcal{F}$ un \'el\'ement de $_{\mathcal{A}_{\vert X}}X_{\mathcal{T}}$, le complexe $R \Gamma(\overline{X}_{\mathcal{T}}, j_{!}\  \mathcal{F})$ est ind\'ependant de la compactification $\overline{X}$ de $X$, et sera not\'e}

$$
R \Gamma_{\textrm{synt},c}(X, \mathcal{F})\ ;
$$

\noindent \textit{ses groupes de cohomologie seront not\'es}

$$
H^i_{\textrm{synt},c}(X, \mathcal{F})
$$

\noindent \textit{et appel\'es groupes de cohomologie syntomique \`a supports compacts.}

\vskip 3mm
\noindent \textit{D\'emonstration}. Si 
$$\xymatrix{
X \ar@{^{(}->}[rr]^{j'}  \ar[rd] & & \overline{X}' \ar[ld] \\
 &  Spec\ k &
 }$$
est une autre compactification de $X$, on note $\overline{X}''$ l'image sch\'ematique de $X$ plong\'e diagonalement dans $\overline{X}\times_{k} \overline{X}'$, $j'' : X \hookrightarrow \overline{X}''$ l'immersion ouverte et $g : \overline{X}'' \rightarrow \overline{X}$, $g' : \overline{X}'' \rightarrow \overline{X}'$ les deux projections (propres).\\

Il s'agit de montrer que\\

\noindent (2.2) $\qquad \qquad \qquad R \Gamma(\overline{X}_{\mathcal{T}}, j_{!}\  \mathcal{F}) = R \Gamma(\overline{X}''_{\mathcal{T}}, j''_{!}\  \mathcal{F}).$\\

Ceci va r\'esulter de la proposition plus g\'en\'erale suivante.

\vskip 3mm
\noindent \textbf{Proposition (2.3)}. 
\textit{Supposons donn\'e un carr\'e cart\'esien de $k$-sch\'emas}

$$
\xymatrix{
X' \ar@{^{(}->}[r]^{j'} \ar[d]_f & \overline{X}' \ar[d]^{\overline{f}}&\\
X \ar@{^{(}->}[r]^{j} & \overline{X}&,
}
$$

\noindent \textit{o\`u $\overline{X}$ est propre sur $k$, $\overline{f}$ est propre, $j$, $j'$ sont des immersions ouvertes. Si $\mathcal{F}$ est un faisceau ab\'elien de torsion sur $X'_{\mathcal{T}}$, alors on a des isomorphismes}\\

$\qquad \qquad R \Gamma(\overline{X}_{\mathcal{T}}, j_{!}\  Rf_{\ast}\ \mathcal{F}) \simeq R \Gamma(\overline{X}_{\mathcal{T}}, R \overline{f}_{\ast}\ j'_{!}\  \mathcal{F})$\\

$\qquad \qquad \qquad \qquad \qquad  \qquad \simeq R \Gamma(\overline{X}'_{\mathcal{T}},  j'_{!}\  \mathcal{F}).$\\

En effet (2.1) r\'esulte de (2.3) via le lemme suivant :

\vskip 3mm
\noindent \textbf{Lemme (2.4)}. 
\textit{Si} 
$$
\xymatrix{
& \overline{X}'' \ar[dd]^{g}\\
X \ar@{^{(}->}[ur]^{j''} \ar@{^{(}->}[rd]_{j}\\
& \overline{X}
}
$$
\textit{est un triangle commutatif de sch\'emas, avec $j$, $j''$ des immersions ouvertes dominantes, alors $X$ est le produit fibr\'e $X \times_{\overline{X}} \overline{X}''$.}

\vskip 3mm
\noindent \textit{D\'emonstration de (2.4)}. Soit $U''$ le produit fibr\'e

$$
\xymatrix{
X \ar@/^/[rrd]^{j''} \ar@/_/[rdd]_{\textrm{id}} \ar@{->}[rd]^{\varphi}\\
& U'' \ar[d]^{f} \ar@{^{(}->}[r]_{\overline{j}} & \overline{X}''  \ar[d]^{g}&\\
& X \ar@{^{(}->}[r]_{j} & \overline{X}&.
}
$$

\noindent Puisque $f \circ \varphi = \textrm{id}$, $\varphi$ est une immersion ferm\'ee [EGA I, (4.3.6) (iv)] ; or $\varphi$ est \'etale, car $\tilde{j} \circ \varphi = j''$ et $\tilde{j}$, $j''$ sont \'etales [EGA IV, (17.3.5)]. Ainsi $\varphi$ est une immersion ouverte [EGA IV, (17.9.1); EGA I, 4.2]. De plus $\varphi$ est dominante car $j''$ l'est ; donc $\varphi$ est surjective, car $\varphi$ est finie. Une immersion ouverte surjective est un isomorphisme. $\square$

\vskip 3mm
\noindent \textit{D\'emonstration de (2.3)}. Faisons la d\'emonstration dans le cas des gros topos syntomiques : pour les petits topos cela r\'esulte de l'\'egalit\'e 

$$
H^i(\overline{X}_{\textrm{SYNT}}, \mathcal{G}) = H^i(\overline{X}_{\textrm{synt}}, \mathcal{G})
$$

\noindent valable pour tout faisceau ab\'elien $\mathcal{G}$ sur synt$(X)$, et tout entier $i \geqslant 0$ [E-LS 2, 1.2] et [Mi, II, prop 3.1]. On d\'esigne par un "ET" en indice les gros topos \'etales [E-LS 2, \S\ 1].\\

On a un cube commutatif de morphismes de gros topos

$$
\xymatrix{
&X'_{\textrm{SYNT}} \ar@{^{(}->}[rr]^{j'_{\textrm{SYNT}}} \ar@{.>}[dd]^(.7){f_{\textrm{SYNT}}} |\hole \ar[dl]_{\beta_{X'}}  && \overline{X}'_{\textrm{SYNT}} \ar[dd]^{\overline{f}_{\textrm{SYNT}}}  \ar[dl]^{\beta_{\overline{X}'}} & & \\
 X'_{\textrm{ET}} \ar@{^{(}->}[rr]^(.7){j'_{\textrm{ET}}} \ar[dd]_{f_{\textrm{ET}}} && \overline{X}'_{\textrm{ET}}  \ar[dd]^(.3){\overline{f}_{\textrm{ET}}}  & &\\
& X_{\textrm{SYNT}} \ar@{^{(}.>}[rr]_(.7){j_{\textrm{SYNT}}} \ar@{.>}[dl]_{\beta{X}} && \overline{X}_{\textrm{SYNT}}   \ar[dl]^{\beta_{\overline{X}}} & & Z_{\textrm{SYNT}} \ar[dl]^{\beta_{Z}}  \ar[ll]_{i_{\textrm{SYNT}}} \\
X_{\textrm{ET}} \ar@{^{(}->}[rr]_{{j_{\textrm{ET}}}}  && \overline{X}_{\textrm{ET}} &&  Z_{\textrm{ET}} \ar[ll]^{i_{\textrm{ET}}}\\
}
$$

\noindent o\`u $Z = \overline{X} \setminus X$, et un isomorphisme

$$
R \Gamma(\overline{X}_{\textrm{SYNT}},\  j_{!}\ R f_{\textrm{SYNT}^{\ast}}\  \mathcal{F}) \simeq R \Gamma(\overline{X}_{\textrm{ET}}, R \beta_{\overline{X}^{\ast}}\  j_{\textrm{SYNT !}}\  R f_{\textrm{SYNT}^{\ast}}\  \mathcal{F}).
$$
\noindent Supposons \'etablie la proposition suivante :

\vskip 3mm
\noindent \textbf{Proposition (2.5)}. 
\textit{Avec les notations pr\'ec\'edentes, annelons $\overline{X}_{\textrm{SYNT}}$ par $\mathcal{A}$ et $\overline{X}_{\textrm{ET}}$ par $\mathcal{B} = \beta_{\overline{X}^{\ast}} \mathcal{A}$. Alors pour tout $\mathcal{H} \in {_{\mathcal{A} \vert X}X_{\textrm{SYNT}}}$ on a un isomorphisme}
$$
j_{\textrm{ET}!}\  R \beta_{X^{\ast}}(\mathcal{H}) \displaystyle \mathop{\rightarrow}^{\sim} R \beta_{\overline{X}^{\ast}}\  j_{\textrm{SYNT !}} (\mathcal{H}).
$$

Alors on a des isomorphismes \\

$R \Gamma(\overline{X}_{\textrm{SYNT}},\  j_{!}\ R f_{\textrm{SYNT}^{\ast}}\  \mathcal{F}) \simeq R \Gamma(\overline{X}_{\textrm{ET}}, j_{\textrm{ET}!}\ R \beta_{X^{\ast}}\  R f_{\textrm{SYNT}^{\ast}}\  \mathcal{F})$\\

$\quad \simeq R \Gamma(\overline{X}_{\textrm{ET}}, j_{\textrm{ET}!}\  R f_{\textrm{ET}^{\ast}}\ R \beta_{X'^{\ast}}\   \mathcal{F})$\\

$\quad \simeq R \Gamma(\overline{X}_{\textrm{ET}}, R \overline{f}_{\textrm{ET}^{\ast}}\ j'_{\textrm{ET}!}\   R \beta_{X'^{\ast}}\   \mathcal{F})$ [SGA 4, XVII, lemme 5.1.6] car $\mathcal{F}$ de torsion.\\

$\quad \simeq R \Gamma(\overline{X}_{\textrm{ET}}, R \overline{f}_{\textrm{ET}^{\ast}}\ R \beta_{\overline{X}'^{\ast}}\   j'_{\textrm{SYNT}!}\  \mathcal{F})$ [Prop (2.5)]\\

$\quad \simeq R \Gamma(\overline{X}_{\textrm{ET}},  R \beta_{\overline{X}^{\ast}}\   R \overline{f}_{\textrm{SYNT}^{\ast}}\ j'_{\textrm{SYNT}!}\  \mathcal{F})$ \\

$\quad \simeq R \Gamma(\overline{X}_{\textrm{SYNT}},  R \overline{f}_{\textrm{SYNT}^{\ast}}\ j'_{\textrm{SYNT}!}\  \mathcal{F})\ ;$ \\

\noindent d'o\`u la proposition (2.3).\\

\noindent\textit{Etablissons la proposition (2.5)}.\\
Puisque $j_{\textrm{ET}!}$ et $j_{\textrm{SYNT}!}$ sont exacts, il suffit de montrer l'isomorphisme \hfill\break $j_{\textrm{ET}!}\  \beta_{X^{\ast}} \displaystyle \mathop{\rightarrow}^{\sim} \beta_{\overline{X}^{\ast}} j_{\textrm{SYNT}!}$.\\
Par le lemme du serpent appliqu\'e au morphisme de suites exactes \\

$$
\begin{array}{c}
\xymatrix{
0 \ar[r] & j_{\textrm{ET}!} \beta_{X^{\ast}} \mathcal{H} \ar[r] \ar@{^{(}->}[ddd] & j_{\textrm{ET}^{\ast}} \beta_{X^{\ast}} \mathcal{H} \ar[r] \ar[ddd]^\simeq]& i_{\textrm{ET}^{\ast}} i^{\ast}_{\textrm{ET}} j_{\textrm{ET}^{\ast}} \beta_{X^{\ast}} \mathcal{H} \ar[r]  \ar[d]^\simeq & 0\\
& & &  i_{\textrm{ET}^{\ast}} i^{\ast}_{\textrm{ET}} \beta_{\overline{X}^{\ast}} j_{\textrm{SYNT}^{\ast}} \mathcal{H} \ar[d] &\\
& & & i_{\textrm{ET}^{\ast}} \beta_{Z^{\ast}} i_{\textrm{SYNT}}^{\ast} j_{\textrm{SYNT}^{\ast}} \mathcal{H} \ar[d]^\simeq &\\
0 \ar[r] & \beta_{\overline{X}^{\ast}} j_{\textrm{SYNT}!} \mathcal{H} \ar[r] & \beta_{\overline{X}^{\ast}} j_{\textrm{SYNT}^{\ast}} \mathcal{H} \ar[r]  & \beta_{\overline{X}^{\ast}} i_{\textrm{SYNT}^{\ast}}  i_{\textrm{SYNT}}^{\ast}    j_{\textrm{SYNT}^{\ast}} \mathcal{H}\ ,  & 
}
\end{array}
\leqno{(2.5.1)}
$$

\noindent il nous suffit de montrer que, pour tout faisceau $\mathcal{G}$, on a un isomorphisme

\noindent (2.5.2) $\qquad \qquad i_{\textrm{ET}^{\ast}}\  i_{\textrm{ET}}^{\ast}\  \beta_{\overline{X}^{\ast}}(\mathcal{G}) \displaystyle \mathop{\rightarrow}^{\varphi \atop{\sim}} i_{\textrm{ET}^{\ast}}\  \beta_{Z^{\ast}}\  i^{\ast}_{\textrm{SYNT}}(\mathcal{G}).$\\

\noindent Or $i_{\textrm{ET}^{\ast}}\ i_{\textrm{ET}}^{\ast}\ Ê\beta_{\overline{X}^{\ast}}(\mathcal{G})$ est le faisceau associ\'e au pr\'efaisceau\\

\vskip 4mm

$\qquad \qquad X' \mapsto i_{\textrm{ET}}^{\ast}\  \beta_{\overline{X}^{\ast}}(\mathcal{G}) (Z \times_{\overline{X}} \overline{X}') = \displaystyle \mathop{\lim}_{\rightarrow \atop{X''}} \beta_{\overline{X}^{\ast}}(\mathcal{G}) (X'')$

$\qquad \qquad  \qquad  \qquad \qquad \qquad \qquad \qquad  \quad = \displaystyle \mathop{\lim}_{\rightarrow \atop{X''}} \mathcal{G} (X'') \simeq \mathcal{G}(Z \times_{\overline{X}} \overline{X}')$\\

\noindent o\`u $\overline{X}'$ est un $\overline{X}$-sch\'ema et la limite inductive est prise sur les diagrammes commutatifs

$$
\xymatrix{
X''  \ar[d] & Z \times_{\overline{X}} \overline{X}' \ar[l]\ \ar[d]\\
\overline{X} & Z \ar[l]
}
$$

\noindent o\`u $X''$ est un $\overline{X}$-sch\'ema ; comme $i_{\textrm{ET}^{\ast}}\ \beta_{Z^{\ast}}\  i^{\ast}_{\textrm{SYNT}}(\mathcal{G})$ a la m\^eme description, il en r\'esulte que $\varphi$ est un isomorphisme. $\square$

\vskip 3mm
\noindent \textbf{Th\'eor\`eme (2.6)}. 
\textit{Soient $i_{1} : Z \hookrightarrow X$ une immersion ferm\'ee entre deux $k$-sch\'emas s\'epar\'es de type fini, $j_{1} : \cup \hookrightarrow X$ l'immersion ouverte du compl\'ementaire et $\mathcal{F} \in {_{\mathcal{A}}X_{\mathcal{T}}}$, o\`u $\mathcal{A}$ est un faisceau d'anneaux de torsion. Alors on a une suite exacte longue de cohomologie syntomique \`a supports compacts }\\

$ \rightarrow H^i_{\textrm{synt},c}(U, \mathcal{F}_{\vert U}) \rightarrow H^i_{\textrm{synt},c}(U, \mathcal{F}) \rightarrow H^i_{\textrm{synt},c}(Z, \mathcal{F}_{\vert Z}) \rightarrow H^{i+1}_{\textrm{synt},c}(U, \mathcal{F}_{\vert U}) \rightarrow $

\vskip 3mm
\noindent \textit{D\'emonstration}. Choisissons une compactification $\overline{X}$ de $X$ au-dessus de $k$, $\overline{j} : X \hookrightarrow \overline{X}$ l'immersion ouverte dominante et soit $\overline{Z}$ l'adh\'erence sch\'ematique de $Z$ dans $\overline{X}$. On a alors un diagramme commutatif \`a carr\'es cart\'esiens

$$
\begin{array}{c}
\xymatrix{
Z \ar@{^{(}->}[r]^{\overline{j}'} \ar@{^{(}->}[d]_{i_{1}} & \overline{Z} \ar@{^{(}->}[d]^{i}\\
X \ar@{^{(}->}[r]_{\overline{j}} & \overline{X}\\
U \ar@{=}[r] \ar@{^{(}->}[u]^{j_{1}} & U \ar@{^{(}->}[u]_{j}
}
\end{array}
\leqno{(2.6.1)}
$$

\noindent o\`u $i$ est une immersion ferm\'ee et $j$, $\overline{j}'$ des immersions ouvertes.\\
En appliquant le foncteur exact $\overline{j}_{!}$ \`a la suite exacte

$$
0 \longrightarrow j_{1 !}\  j^{\ast}_{1}\  \mathcal{F} \longrightarrow \mathcal{F} \longrightarrow\  i_{1^{\ast}}\  i^{\ast}_{1}\  \mathcal{F} \longrightarrow 0,
$$

\noindent on obtient la suite exacte\\

\noindent (2.6.2) $\qquad \qquad 0 \longrightarrow j_{!}\  j^{\ast}_{1}\   \mathcal{F} \longrightarrow \overline{j}_{!}\ \mathcal{F} \longrightarrow  \overline{j}_{!}\  i_{1^{\ast}}\  i^{\ast}_{1}\  \mathcal{F} \longrightarrow 0$\\

\noindent car $\overline{j}_{!}\ j_{1 !} = j_{!}$ [Th\'eo (1.5) (1)]. Or l'exactitude des foncteurs $i_{1 \ast}$ et $i_{\ast}$, jointe \`a la proposition (2.3), donne un isomorphisme

$$
R \Gamma(\overline{X}_{\mathcal{T}}, \overline{j}_{!}\ i_{1 \ast}\ i^{\ast}_{1}\ \mathcal{F} \displaystyle \mathop{\longrightarrow}^{\sim} R \Gamma (\overline{Z}_{\mathcal{T}}, \overline{j'}_{!}\ i^{\ast}_{1}\ \mathcal{F})\ ;
$$

\noindent ainsi, par application du foncteur $R \Gamma(\overline{X}_{\mathcal{T}}, -)$ \`a la suite exacte (2.6.2) on obtient, via (2.1), un triangle distingu\'e\\

\noindent (2.6.3) $\qquad R \Gamma_{\textrm{synt},c}(U, \mathcal{F}_{\vert U}) \longrightarrow R \Gamma_{\textrm{synt},c}(X, \mathcal{F}) \longrightarrow R \Gamma_{\textrm{synt},c}(Z, \mathcal{F}_{\vert Z}),$\\

\noindent qui fournit \`a son tour la suite exacte longue du th\'eor\`eme. $\square$\\

\vskip 3mm
\section*{3. Comparaison avec la cohomologie \'etale et la cohomologie rigide}

\textbf{3.0.} On suppose dans ce $\S$\  3 que le corps $k$ contient $\mathbb{F}_{q},\ q = p^a$. On d\'esigne par $C(k)$ un anneau de Cohen de $k$ de caract\'eristique 0 et par $K_{0}$ le corps des fractions de $C(k)$.\\

Comme en [Et 3, III (3.3)] ou [Et 5, (3.3)], $\mathcal{V}$ est un anneau de valuation discr\`ete complet, d'uniformisante $\pi$ et $\sigma : \mathcal{V \rightarrow \mathcal{V}}$ un rel\`evement de la puissance $q$ de $k$ tel que $\sigma(\pi) = \pi$ construit via [Et 2, 1.1].\\

On note $e$ l'indice de ramification de $\mathcal{V}$, $K = \textrm{Frac}(\mathcal{V})$, $\mathcal{V}^{\sigma} = \textrm{Ker} \{1 - \sigma : \mathcal{V} \rightarrow \mathcal{V} \}, \mathcal{V}_{n} = \mathcal{V} / \pi^{n+1} \mathcal{V}, \mathcal{V}^{\sigma}_{n} = \mathcal{V}^{\sigma} / \pi^{n+1} \mathcal{V}^{\sigma}$ et $K^{\sigma} = \textrm{Frac} (\mathcal{V}^{\sigma})$.\\

Si $X$ est un sch\'ema on dit qu'un $\mathcal{V}^{\sigma}_{n}$-module $\mathcal{F}$ sur \'et$(X)$ est localement trivial (on dit aussi constant-tordu constructible ou encore localement constant constructible) s'il est localement isomorphe \`a une somme directe finie de copies de $\mathcal{V}^{\sigma}_{n}$ : c'est alors la m\^eme chose de dire qu'il est localement trivial sur SYNT$(X)$ [E-LS 2, 5.1].\\

Un $\mathcal{V}^{\sigma}$-faisceau lisse $\mathcal{F}$ (localement libre de rang fini) sur $X$ est un syst\`eme projectif $\mathcal{F} = (\mathcal{F}_{n})_{n \in \mathbb{N}}$ o\`u, pour tout $n$, $\mathcal{F}_{n}$ est un $\mathcal{V}^{\sigma}_{n}$-module localement trivial sur \'et$(X)$ et pour $n' \geqslant n$, $\mathcal{F}_{n} = \mathcal{F}_{n'} \otimes \mathcal{V}^{\sigma}_{n}$. Les $K^{\sigma}$-faisceaux lisses sont les $\mathcal{V}^{\sigma}$-faisceaux lisses \`a isog\'enie pr\`es.\\

Pour un $\mathcal{V}^{\sigma}$-faisceau lisse $\mathcal{F}$ (resp. un $K^{\sigma}$-faisceau lisse   $\mathcal{F_{\mathbb{Q}}}$) sur $X$ on pose, pour tout entier $i \geqslant 0$

$$
H^i_{\textrm{synt},c}(X, \mathcal{F}) := \displaystyle \mathop{\lim}_{\leftarrow \atop_{n}} H^i_{\textrm{synt},c}(X, \mathcal{F}_{n})
$$\\

[resp. $ \qquad \qquad  H^i_{\textrm{synt},c}(X, \mathcal{F}_{\mathbb{Q}}) = (\displaystyle \mathop{\lim}_{\leftarrow \atop_{n}} H^i_{\textrm{synt},c}(X, \mathcal{F}_{n})) \otimes_{\mathcal{V}^{\sigma}} K^{\sigma}]\ .
$

\vskip 3mm
\noindent \textbf{Proposition (3.1)}. 
\textit{Soient $X$ un sch\'ema et $\mathcal{F}$ un faisceau ab\'elien sur \'et$(X)$ ; on note encore $\mathcal{F}$ son image inverse par le morphisme de topos $X_{\textrm{synt}} \rightarrow X_{\textrm{\'et}}$. Alors, pour tout entier $i \geqslant 0$, on a un isomorphisme}

$$
H^i_{\textrm{\'et},c} (X, \mathcal{F}) \displaystyle \mathop{\longrightarrow}^{\sim} H^i_{\textrm{synt},c}(X, \mathcal{F}).
$$

\vskip 3mm
\noindent \textit{D\'emonstration}. R\'esulte de [E-LS 2, 1.3]. \ $\square$

\vskip 3mm
\noindent \textbf{Th\'eor\`eme (3.2)}. 
\textit{Supposons $k$ s\'eparablement clos. Soient $X$ un $k$-sch\'ema s\'epar\'e de type fini, $\mathcal{F}_{\mathbb{Q}}$ un $K^{\sigma}$-faisceau lisse sur $X$ et $f : Y \rightarrow X$ sur $k$-morphisme fini \'etale galoisien de groupe $G$. Alors, pour tout entier $i \geqslant 0$, on a des isomorphismes de $K^{\sigma}$-espaces vectoriels de dimension finie}

$$
H^i_{\textrm{\'et},c}(X, \mathcal{F}_{\mathbb{Q}}) \simeq H^i_{\textrm{\'et},c}(X, f_{\ast} f^{\ast} \mathcal{F}_{\mathbb{Q}})^G \simeq H^i_{\textrm{\'et},c}(Y, f^{\ast} \mathcal{F}_{\mathbb{Q}})^G
$$

$$
\simeq H^i_{\textrm{synt},c}(X, \mathcal{F}_{\mathbb{Q}}) \simeq H^i_{\textrm{synt},c} (X, f_{\ast} f^{\ast} \mathcal{F}_{\mathbb{Q}})^G \simeq H^i_{\textrm{synt},c}(Y, f^{\ast} \mathcal{F}_{\mathbb{Q}})^G.
$$

\vskip 3mm
\noindent \textit{D\'emonstration}. Compte tenu de (3.1) il suffit de montrer l'assertion pour la cohomologie \'etale. Puisque $f$ est fini, $f_{\ast}$ est exact, et on est ramen\'e \`a montrer l'isomorphisme\\

\noindent (3.2.1) $\qquad  H^i_{\textrm{\'et},c}(X, \mathcal{F}) \otimes_{\mathcal{V}^{\sigma}} K^{\sigma} \simeq [H^i_{\textrm{\'et},c}(X, f_{\ast} f^{\ast} \mathcal{F}) \otimes_{\mathcal{V}^{\sigma}} K^{\sigma}]^G $ \\

\noindent pour un $\mathcal{V}^{\sigma}$-faisceau lisse $\mathcal{F}$.\\
Comme $\mathcal{F}$ est localement libre on \'etablit le lemme suivant comme [Et 1, III (3.1.2)].

\vskip 3mm
\noindent \textbf{Lemme (3.2.2)}. 
\textit{Sous les hypoth\`eses pr\'ec\'edentes, on a un isomorphisme}

$$
\mathcal{F} \displaystyle \mathop{\longrightarrow}^{\sim} (f_{\ast} f^{\ast}(\mathcal{F}))^G.
$$    

Soient $\overline{X}$ une compactification de $X$ au-dessus de $k$ et $j : X \hookrightarrow \overline{X}$ l'immersion ouverte correspondante. Par exactitude du foncteur $j_{!}$ on d\'eduit de (3.2.2) des isomorphismes\\

\noindent (3.2.3) $\qquad \qquad  j_{!}   \mathcal{F}   \displaystyle \mathop{\longrightarrow}^{\sim} j_{!} (f_{\ast} f^{\ast} \mathcal{F})^G \simeq (j_{!} f_{\ast} f^{\ast} \mathcal{F})^G.$ \\

\noindent Le corps $k$ \'etant s\'eparablement clos, les groupes $H^i_{\textrm{\'et},c}(X, \mathcal{F}_{n})$ et   $H^i_{\textrm{\'et},c}(X, f_{\ast} f^{\ast}  \mathcal{F}_{n})$ sont des $\mathcal{V}^{\sigma}_{n}$-modules de type fini [SGA 4, XVII, 5.3.8] ; par suite les groupes $H^i_{\textrm{\'et},c}(X, \mathcal{F})$ et    $H^i_{\textrm{\'et},c}(X,  f_{\ast} f^{\ast}  \mathcal{F})$ sont des $\mathcal{V}^{\sigma}$-modules de type fini, engendr\'es par tout sous-ensemble qui les engendre mod.$\pi$. \\

On ach\`eve la d\'emonstration de (3.2) comme [Et 1, III, 3.1.1]. $\square$  

\vskip 3mm
\textbf{3.3.} Soient $X$ un $k$-sch\'ema et $\mathcal{H}$ un $\mathcal{V}^{\sigma}_{n}$-module localement trivial sur SYNT$(X)$. On consid\`ere le morphisme de topos annel\'es [E-LS 2, 5.3]

$$
u = u^{(m)}_{X/\mathcal{V}_{n}-\textrm{SYNT}} : ((X/\mathcal{V}_{n})^{(m)}_{\textrm{CRIS-SYNT}}, \mathcal{O}^{(m)}_{X/\mathcal{V}_{n}}) \longrightarrow (X_{\textrm{SYNT}}, \mathcal{V}^{\sigma}_{n})
$$

\noindent et on note

$$
T^{(m)}(\mathcal{H}) := u^{\ast}(\mathcal{H}) = \mathcal{H} \otimes_{\mathcal{V}^{\sigma}_{n}} \mathcal{O}^{(m)}_{X/\mathcal{V}_{n}},
$$

\noindent et
$$T^{(m)}(\mathcal{H})^{\textrm{cris}} := u_{\ast} (T^{(m)}(\mathcal{H})) = u_{\ast} u^{\ast}(\mathcal{H})\ [\mbox{E-LS 2}, 1.11]$$

\noindent (3.3.1) $\qquad \qquad \qquad = \mathcal{H}  \otimes_{\mathcal{V}^{\sigma}_{n}} u_{\ast}({\mathcal{O}}^{(m)}_{X/\mathcal{V}_{n}}) = \mathcal{H}  \otimes_{\mathcal{V}^{\sigma}_{n}} \mathcal{O}^{m-\textrm{cris}}_{n, X},$\\

\noindent o\`u l'on a pos\'e $\mathcal{O}^{m-\textrm{cris}}_{n, X} := u_{\ast}(\mathcal{O}^{(m)}_{X/\mathcal{V}_{n}}).$\\

Soient $i : Z \hookrightarrow X$  une immersion  ferm\'ee et $j : U \hookrightarrow X$ l'immersion ouverte du compl\'ementaire : $j$ et $i$ d\'efinissent respectivement des morphismes de topos [cf (1.1)]

$$
j : U_{\textrm{SYNT}} \longrightarrow X_{\textrm{SYNT}}\ ,
$$

$$
i : Z_{\textrm{SYNT}} \longrightarrow X_{\textrm{SYNT}}\ ,
$$

\noindent et m\^eme des morphismes de topos annel\'es\\

\noindent (3.3.2) $\qquad \qquad \qquad j_{\mathcal{V}} : (U_{\textrm{SYNT}}, \mathcal{V}^{\sigma}_{n,U}) \longrightarrow (X_{\textrm{SYNT}}, \mathcal{V}^{\sigma}_{n,X})$\\

\noindent (3.3.3) $\qquad \qquad \qquad i_{\mathcal{V}} : (Z_{\textrm{SYNT}}, \mathcal{V}^{\sigma}_{n,Z}) \longrightarrow (X_{\textrm{SYNT}}, \mathcal{V}^{\sigma}_{n,X})$\\

\noindent o\`u $\mathcal{V}^{\sigma}_{n,U} = j^{-1}(\mathcal{V}^{\sigma}_{n,X}) \quad, \quad \mathcal{V}^{\sigma}_{n,Z} = i^{-1}(\mathcal{V}^{\sigma}_{n,X}),$\\

\noindent $\mathcal{V}^{\sigma}_{n,X}$ \'etant le faisceau d'anneaux $\mathcal{V}^{\sigma}_{n}$ sur $X_{\textrm{SYNT}}.$\:

On d\'eduit alors de (1.4) l'existence de six foncteurs

$$
\begin{array}{c}
\xymatrix{
& \ar[l]_{i^{\ast}_{\mathcal{V}}} \qquad & \qquad  \ar[l]_{j_{\mathcal{V}!}} &\\
\ _{\mathcal{V}^{\sigma}_{n,Z} }Z_{\textrm{SYNT}} \ar[r]^{i_{\mathcal{V}^{\ast}}} & \ _{\mathcal{V}^{\sigma}_{n,X} }X_{\textrm{SYNT}} \ar[r]^{j^{\ast}_{\mathcal{V}}} & \  _{\mathcal{V}^{\sigma}_{n,U} }U_{\textrm{SYNT}}\\
 &\  \ar[l]_{i^{!}_{\mathcal{V}}} \  &  \  \ar[l]_{j_{\mathcal{V}^{\ast}}}  \ &
}
\end{array}
\leqno{(3.3.4)}
$$

\vskip 3mm
\noindent \textbf{Lemme (3.3.5)}. 
\textit{Sous les hypoth\`eses pr\'ec\'edentes on a des isomorphismes de faisceaux sur les gros sites syntomiques :}\\

\noindent (3.3.5.1) $\qquad \qquad j^{-1} \mathcal{O}^{m-\textrm{cris}}_{n,X} \displaystyle \mathop{\longrightarrow}^{\sim} \mathcal{O}^{m-\textrm{cris}}_{n,U}$\\

\noindent (3.3.5.2) $\qquad \qquad i^{-1} \mathcal{O}^{m-\textrm{cris}}_{n,X} \displaystyle \mathop{\longrightarrow}^{\sim} \mathcal{O}^{m-\textrm{cris}}_{n,Z}.$\\

\vskip 3mm
\noindent \textit{D\'emonstration}. Comme $j^{-1} \mathcal{O}^{m-\textrm{cris}}_{n,X}$ est le faisceau associ\'e au pr\'efaisceau qui \`a tout $U$-sch\'ema $U'$ associe\\
$$
\mathcal{O}^{m-\textrm{cris}}_{n,X}(U') = \Gamma((U'/\mathcal{V}_{n})^{(m)}_{\textrm{CRIS-SYNT}}, \mathcal{O}^{(m)}_{U'/\mathcal{V}n})\  [\mbox{E-LS 2}, 1.10] 
$$
$\qquad \qquad \qquad \qquad \qquad  = \mathcal{O}^{m-\textrm{cris}}_{n,U}(U'),$\\

\noindent on a bien (3.3.5.1). De m\^eme pour (3.3.5.2). $\square$\\

Gr\^ace au lemme (3.3.5) les morphismes de topos $j$ et $i$ ci-dessus induisent des morphismes de topos annel\'es\\

\noindent (3.3.6) $\qquad \qquad j_{\mathcal{O}} : (U_{\textrm{SYNT}}, \mathcal{O}^{m-\textrm{cris}}_{n,U}) \longrightarrow (X_{\textrm{SYNT}}, \mathcal{O}^{m-\textrm{cris}}_{n,X})$\\

\noindent (3.3.7) $\qquad \qquad i_{\mathcal{O}} : (Z_{\textrm{SYNT}}, \mathcal{O}^{m-\textrm{cris}}_{n,Z}) \longrightarrow (X_{\textrm{SYNT}}, \mathcal{O}^{m-\textrm{cris}}_{n,X})$\\

\noindent et six foncteurs

$$
\begin{array}{c}
\xymatrix{
& \ar[l]_{i^{\ast}_{\mathcal{O}}} \qquad & \qquad  \ar[l]_{j_{\mathcal{O}!}} &  {} &\\
\  _{\mathcal{O}^{m-\textrm{cris}}_{n,Z} } Z_{\textrm{SYNT}} \ar[r]^{i_{\mathcal{O}^{\ast}}} & \  _{\mathcal{O}^{m-\textrm{cris}}_{n,X} } X_{\textrm{SYNT}} \ar[r]^{j^{\ast}_{\mathcal{O}}} & \  _{\mathcal{O}^{m-\textrm{cris}}_{n,U} } U_{\textrm{SYNT}}&\\
 &\  \ar[l]_{i^{!}_{\mathcal{O}}} \ &\   \ar[l]_{j_{\mathcal{O}^{\ast}}}\ & .
}
\end{array}
\leqno{(3.3.8)}
$$
    
En utilisant la description (1.4), ou le th\'eor\`eme (1.5) \textit{(6.1)} on en d\'eduit la proposition suivante :    

\vskip 3mm
\noindent \textbf{Proposition (3.3.9)}. 
\textit{Sous les hypoth\`eses et notations de (3.3) on a :}
\begin{itemize}
\item[(1)] \textit{Si $\mathcal{G}$ est un $\mathcal{O}^{m-\textrm{cris}}_{n,U}$-module, alors on a un isomorphisme canonique}\\
$$j_{\mathcal{V}!}(\mathcal{G}) \displaystyle \mathop{\rightarrow}^{\sim} j_{\mathcal{O}!}(\mathcal{G}).$$
\item[(2)] \textit{Si $\mathcal{G}$ est un $\mathcal{V}^{\sigma}_{n,U}$-module, alors on a des isomorphismes canoniques}
$$j_{\mathcal{V}!}(\mathcal{G} \otimes_{\mathcal{V}^{\sigma}_{n,U}} \mathcal{O}^{m-\textrm{cris}}_{n,U}) \simeq j_{\mathcal{O}!} (\mathcal{G} \otimes_{\mathcal{V}^{\sigma}_{n,U}} \mathcal{O}^{m-\textrm{cris}}_{n,U})$$
$\qquad \qquad \simeq j_{\mathcal{V}!}(\mathcal{G}) \otimes_{\mathcal{V}^{\sigma}_{n,X}} j_{\mathcal{V}!} (\mathcal{O}^{m-\textrm{cris}}_{n,U}) \simeq j_{\mathcal{V}!} (\mathcal{G}) \otimes_{\mathcal{V}^{\sigma}_{n,X}} j_{\mathcal{O}!} (\mathcal{O}^{m-\textrm{cris}}_{n,U})$\\

$\qquad \qquad \simeq  j_{\mathcal{V}!}(\mathcal{G}) \otimes_{\mathcal{V}^{\sigma}_{n,X}} \mathcal{O}^{m-\textrm{cris}}_{n,X}.$

\end{itemize} 

\vskip 3mm
Les formules du (2) sont \`a comparer \`a celles de [SGA 4, IV, prop 12.11 (b)]. \\

Supposons \`a pr\'esent que $\mathcal{F}$ est un $\mathcal{V}^{\sigma}_{n}$-module localement trivial sur SYNT$(U) : \mathcal{F}$ \'etant localement trivial, le morphisme $F^{\ast} : \mathcal{F} \rightarrow \mathcal{F}^{(q)} = F^{-1}_{X}(\mathcal{F})$ est un isomorphisme ; on pose $F = F^{\ast-1} : \mathcal{F}^{(q)} \displaystyle \mathop{\rightarrow}^{\sim} \mathcal{F}$ et $\phi _{U}= T^{(m)}(F): T^{(m)}(\mathcal{F})^{(q)}=T^{(m)}(\mathcal{F}^{(q)})\displaystyle \mathop{\rightarrow}^{\sim} T^{(m)}(\mathcal{F}) $ qui munit $T^{(m)}(\mathcal{F})$ d'une structure de $F\mbox{-}m$-cristal localement trivial [E-LS 2, 5.3]. \\
\noindent

On note encore $\phi_{U} : T^{(m)}(\mathcal{F})^{\textrm{cris}} \rightarrow T^{(m)}(\mathcal{F})^{\textrm{cris}}$ l'homomorphisme obtenu en composant $\phi_{U}^{\textrm{cris}}$ avec  $F^{\ast} :  T^{(m)}(\mathcal{F})^{\textrm{cris}}   \rightarrow T^{(m)}(\mathcal{F})^{\textrm{cris}(q)}$ [E-LS 2, 5.2]. Alors la suite exacte de [E-LS 2, th\'eo 5.5] s'interpr\`ete, via (3.3.1), comme une suite exacte sur SYNT$(U)$ :  \\

\noindent (3.3.10) $\qquad \qquad 0 \longrightarrow \mathcal{F} \longrightarrow u_{\ast}\  u^{\ast}(\mathcal{F}) \displaystyle \mathop{\longrightarrow} _{1-\phi_{U}} u_{\ast}\ u^{\ast}(\mathcal{F}) \longrightarrow 0$ \\

\noindent ou encore\\

\noindent (3.3.11) $\qquad \quad 0 \longrightarrow \mathcal{F} \longrightarrow \mathcal{F} \otimes_{\mathcal{V}^{\sigma}_{n}}\  \mathcal{O}^{m-\textrm{cris}}_{n,U} \displaystyle \mathop{\longrightarrow} _{1-\phi_{U}}  \mathcal{F} \otimes_{\mathcal{V}^{\sigma}_{n}}\  \mathcal{O}^{m-\textrm{cris}}_{n,U} \longrightarrow 0$\ .\\

\noindent En lui appliquant le foncteur exact $j_{\mathcal{V}!}$ [th\'eo 1.5], on obtient encore une suite exacte, ce qui, compte tenu de (3.3.9), \'etablit le th\'eor\`eme suivant :

\vskip 3mm
\noindent \textbf{Th\'eor\`eme (3.3.12)}. 
\textit{Sous les notations de (3.3), si $\mathcal{F}$ est un $\mathcal{V}^{\sigma}_{n}$-module localement trivial sur SYNT$(U)$, alors on a des suites exactes de $\mathcal{V}^{\sigma}_{n}$-modules sur SYNT$(X)$ :}

$$
\xymatrix{
0 \ar[r] & j_{\mathcal{V}!}\ \mathcal{F} \ar[r] \ar@{=}[d]& j_{\mathcal{O}!}(\mathcal{F} \otimes_{\mathcal{V}^{\sigma}_{n}} \mathcal{O}^{m-\textrm{cris}}_{n,U}) \ar[r]_{1-\phi} \ar[d]^{\simeq} & j_{\mathcal{O}!} (\mathcal{F} \otimes_{\mathcal{V}^{\sigma}_{n}} \mathcal{O}^{m-\textrm{cris}}_{n,U})\ar[r] &  0\\
0 \ar[r] & j_{\mathcal{V}!}\ \mathcal{F} \ar[r] & j_{\mathcal{V}!}(\mathcal{F}) \otimes_{\mathcal{V}^{\sigma}_{n}} \mathcal{O}^{m-\textrm{cris}}_{n,X} \ar[r]_{1-\phi} & j_{\mathcal{V}!} (\mathcal{F}) \otimes_{\mathcal{V}^{\sigma}_{n}} \mathcal{O}^{m-\textrm{cris}}_{n,X} \ar[r] &  0
}
$$

\noindent \textit{o\`u $\phi = j_{\mathcal{V}!} (\phi_{U})$.}\\

Les suites exactes du th\'eor\`eme (3.3.12) vont nous permettre en passant \`a la cohomologie dans le th\'eor\`eme suivant de relier cohomologie \'etale et cohomologie rigide.\\

\noindent \textbf{Th\'eor\`eme (3.3.13)}. 
\textit{On suppose le corps $k$ parfait. Soient $X$ un $k$-sch\'ema s\'epar\'e de type fini, $\mathcal{F}_{\mathbb{Q}}$ un $K^{\sigma}$-faisceau lisse sur $X$ et $E_{K} \in F^a\mbox{-}\textrm{Isoc}(X/K)^{\circ}$ le $F$-isocristal convergent associ\'e \`a $\mathcal{F}_{\mathbb{Q}}$ [E-LS 2 ; 5.6] et on suppose que $E_{K}$ provient de $E^{\dag}_{K} \in F^a\mbox{-}\textrm{Isoc}^{\dag}(X/K)^{\circ}$ par le foncteur d'oubli $F^a\mbox{-}\textrm{Isoc}^{\dag}(X/K)^{\circ} \rightarrow F^a\mbox{-}\textrm{Isoc}(X/K)^{\circ}$. On note $\mathcal{F} = (\mathcal{F}_{n})_{n \in \mathbb{N}}$ un $\mathcal{V}^{\sigma}$-faisceau lisse associ\'e \`a $\mathcal{F}_{\mathbb{Q}}$ et $E^{m-\textrm{cris}}_{n} = T^{(m)}(\mathcal{F}_{n})^{\textrm{cris}}$ [cf (3.3.1)].\\
Alors on a :}
\begin{itemize}
\item[(1)] \textit{Il existe un isomorphisme canonique}
$$
R \Gamma_{\textrm{rig},c}(X/K, E^{\dag}_{K}) = R \displaystyle \mathop{\lim}_{\leftarrow \atop{m}} [(R  \displaystyle \mathop{\lim}_{\leftarrow \atop{n}} R\Gamma_{\textrm{synt},c}(X, E_{n}^{m-\textrm{cris}})) \otimes \mathbb{Q}].
$$
\item[(2)] \textit{Si de plus $k$ est alg\'ebriquement clos, il existe, pour tout entier $i \geqslant 0$, une suite exacte courte}
$$
O \rightarrow H^i_{\textrm{\'et},c}(X, \mathcal{F}_{\mathbb{Q}}) \rightarrow H^i_{\textrm{rig},c}(X/K,E^{\dag}_{K}) \displaystyle \mathop{\longrightarrow}_{1-\phi} H^i_{\textrm{rig},c}(X/K, E^{\dag}_{K}) \rightarrow 0.$$
\end{itemize}

\vskip 3mm
\noindent \textit{D\'emonstration}.\\
\textit{Prouvons le (1)}. Comme la cohomologie rigide \`a supports compacts ne d\'epend que du sch\'ema r\'eduit sous-jacent \`a $X$, on peut supposer $X$ r\'eduit. On va faire une r\'ecurrence sur la dimension de $X$.\\

Si dim $X = 0$, alors $X = \displaystyle \mathop{\bigcup}_{\textrm{finie}} \textrm{Spec} A_{i}$ o\`u $A_{i}$ est artinien [Eis., cor 9.1] car $X$ est de type fini sur $k$, et $A_{i}$ est un produit fini $\displaystyle \mathop \Pi_{j} A_{i,m_{j}}$ d'anneaux artiniens locaux r\'eduits ($X$ est r\'eduit) : ainsi $A_{i,m_{j}}$ est un corps [Bour, A VIII, $\S$\ 6, \no4, prop 9] $k_{ij}$ extension finie du corps parfait $k$ [Eis, cor 2.15]. En particulier $X$ est fini \'etale sur $k$, et alors l'assertion du th\'eor\`eme est prouv\'ee dans [E-LS 2, prop 3.11].\\

Si dim $X \geqslant 1$, il existe, puisque $k$ est parfait et $X$ r\'eduit, un ouvert non vide $U \hookrightarrow X$ qui est lisse sur $k$, de ferm\'e compl\'ementaire $Z$ tel que dim $Z < \textrm{dim}\ X$. Comme les deux foncteurs $R \Gamma_{\textrm{rig},c}(X/K, -)$ et $R \displaystyle \mathop{\lim}_{\leftarrow \atop{m}} \ [(R \displaystyle \mathop{\lim}_{\leftarrow \atop{n}} R\Gamma_{\textrm{synt},c}(X, (-)_{n}^{m-\textrm{cris}})) \otimes \mathbb{Q})]$ donnent lieu \`a des triangles distingu\'es faisant intervenir $X$, $U$ et $Z$ on est ramen\'e \`a prouver le th\'eor\`eme pour $U$. Donc on peut supposer $X$ lisse connexe et aussi $K' = K$ avec $e \leqslant p-1$, avec $\mathcal{F}$ un $\mathcal{V}^{\sigma}$-faisceau lisse sur $X$, associ\'e \`a un $F$-cristal unit\'e $E$ sur $X/\mathcal{V}$, d'isocristal convergent unit\'e $E_{K}$ par [B, (2.4.2)] suppos\'e provenir de $E^{\dag}_{K} \in F^a\mbox{-}\textrm{Isoc}(X/K)^{\circ}$.\\

Notons $\overline{X}$ une compactification de $X$ sur $k$. D'apr\`es le th\'eor\`eme de monodromie finie ``g\'en\'erique'' de Tsuzuki [Tsu 2, theo 3.1] il existe un $k$-sch\'ema projectif et lisse $\overline{X}'$, un $k$-morphisme propre surjectif $\overline{w} : \overline{X'} \rightarrow \overline{X}$ g\'en\'eriquement \'etale, tel qu'en posant $X' = \overline{w}^{-1}(X)$, $j' : X' \hookrightarrow \overline{X'}$ l'immersion ouverte, il existe un unique $N^{\dag} \in F^a\mbox{-}\textrm{Isoc}^{\dag}(\overline{X'}/K)^{\circ}$ avec $\overline{w}^{\ast}(E^{\dag}_{K}) \simeq (j')^{\dag}(N^{\dag})$.\\

\noindent Notons $U \hookrightarrow X$ un ouvert dense tel que la restriction $w : U' = U \times_{\overline{X}} \overline{X'} \rightarrow U$ de $\overline{w}$ soit finie \'etale : quitte \`a r\'etr\'ecir $U$ on peut supposer $U$ affine et int\'egralement clos, de m\^eme pour $U'$. Puisque $U$ et $U'$ sont connexes il existe un morphisme fini \'etale $s : U'' \rightarrow U'$ tel que le compos\'e $f : w \circ s : U'' \rightarrow U$ soit fini \'etale galoisien de groupe not\'e $G$ [Mi, I, Rk 5.4]. D\'esignons par $\overline{X''}$ la fermeture int\'egrale de $\overline{X'}$ dans $U''$ : on obtient un diagramme commutatif \`a carr\'es cart\'esiens

$$
\xymatrix{
U'' \ar@{^{(}->}[rr]^{j''} \ar@/_2pc/[dd]_{f} \ar[d]_{s}  &  & \overline{X''} \ar[d]^{\overline{s}}\\
U' \ar@{^{(}->}[r]  \ar[d]_{w} & X' \ar@{^{(}->}[r] \ar[d] & \overline{X'} \ar[d]^{\overline{w}}\\
U \ar@{^{(}->}[r]  \ar@/_1pc/[rr]_{j_{U}} & X \ar@{^{(}->}[r] ^{j} & \overline{X}
}
$$

\noindent o\`u $\overline{s}$ est un morphisme fini et les fl\`eches horizontales sont des immersions ouvertes : en particulier $\overline{X''}$ est une compactification de $U''$.
On note $\overline{Z}$ l'adh\'erence sch\'ematique de $Z$ dans $\overline{X}$ et $j_{Z} : Z \hookrightarrow \overline{Z}$ l'immersion ouverte dominante.\\

Puisque $\overline{X'}$ est propre et lisse sur $k$, il existe, d'apr\`es (3.3.13) un $F$-cristal unit\'e $M$ sur $\overline{X'}$ tel que $N^{\dag} \simeq M^{an}$. Posons $\overline{\mathcal{M}} = \overline{s}^{\ast}_{\textrm{CRIS}}(M)$, $\mathcal{M} = j''^{\ast}_{\textrm{CRIS}}(\overline{\mathcal{M}})$ ; d'apr\`es [B, (2.42)] $\mathcal{M}$ est isog\`ene \`a $f^{\ast}_{\textrm{CRIS}}(E_{\vert U})$, i.e. il existe un entier $r \geqslant 0$ et des morphismes

$$
\alpha : \mathcal{M} \longrightarrow f^{\ast}_{\textrm{CRIS}}(E_{\vert U}),
$$

$$
\beta : f^{\ast}_{\textrm{CRIS}}(E_{\vert U}) \longrightarrow \mathcal{M},
$$

\noindent tels que $\alpha \circ \beta = p^r$ et $\beta \circ \alpha = p^r$. On notera $\overline{\mathcal{M}}_{K}$ le $F$-isocristal (sur)convergent sur $\overline{X''}$ associ\'e \`a $\overline{\mathcal{M}}$ et $\mathcal{M}^{\dag}_{K} = j''^{\dag}(\overline{\mathcal{M}}_{K})$. D'apr\`es [B-M 1, cor du th\'eo 6] le $F$-cristal unit\'e $\overline{\mathcal{M}}$ est le cristal de Dieudonn\'e d'un groupe $p$-divisible \'etale, dont le $\mathcal{V}^{\sigma}$-faisceau lisse associ\'e sera not\'e $\overline{\mathcal{G}} = (\overline{\mathcal{G}}_{n})_{n \in \mathbb{N}}$ ; on pose $\mathcal{G} = j''^{\ast}(\overline{\mathcal{G}})$.

L'isog\'enie $\alpha$ (resp $\beta$) fournit une isog\'enie $\alpha_{\mathcal{G}} : \mathcal{G} \rightarrow f^{\ast}(\mathcal{F}_{\vert U})$ (resp $\beta_{\mathcal{G}} : f^{\ast}(\mathcal{F}_{\vert U}) \rightarrow \mathcal{G})$ telle que $\alpha_{\mathcal{G}} \circ \beta_{\mathcal{G}} = p^r$ et $\beta_{\mathcal{G}} \circ \alpha_{\mathcal{G}} = p^r$.

L'isog\'enie $\alpha$ (resp $\beta$) fournit, par la construction de Berthelot [B, (2.4.2)], un isomorphisme sur les $F$-isocristaux convergents associ\'es

$$
\alpha_{K} : \mathcal{M}_{K} \displaystyle \mathop{\longrightarrow}^{\sim} f^{\ast}_{\textrm{rig}}(E_{K \vert U})
$$

$$
(\textrm{resp}\  \beta_{K} : f^{\ast}_{\textrm{rig}}(E_{K \vert U}) \displaystyle \mathop{\longrightarrow}^{\sim} \mathcal{M}_{K}).
$$
\noindent D'apr\`es [Et 2, th\'eo 5] l'isomorphisme $\alpha_{K}$ (resp. $\beta_{K}$) se rel\`eve de mani\`ere unique en un isomorphisme
$$
\alpha_{K}^{\dag} : \mathcal{M}_{K}^{\dag} \displaystyle \mathop{\longrightarrow}^{\sim} f^{\ast}_{\textrm{rig}}(E_{K \vert U}^{\dag})
$$

$$
(\textrm{resp}\  \beta_{K}^{\dag} : f^{\ast}_{\textrm{rig}}(E_{K \vert U}^{\dag}) \displaystyle \mathop{\longrightarrow}^{\sim} \mathcal{M}_{K}^{\dag});
$$

\noindent de m\^eme l'action de $G$ sur $f^{\ast}_{\textrm{rig}}(E_{K \vert U})$ se rel\`eve de mani\`ere unique \`a $f^{\ast}_{\textrm{rig}}(E_{K \vert U}^{\dag})$.\\

D'autre part, par le th\'eor\`eme (3.3.12), on a un morphisme de suites exactes sur SYNT$(\overline{X''})$

$$
\xymatrix{
 (S_{1})\ 0 \ar[r]  & j_{!} f^{\ast}(\mathcal{F}_{n \vert U}) \ar[r] \ar[d]_{{j''_{!}}(\beta_{\mathcal{G}})_{n}} & j_{!} f^{\ast}_{\textrm{CRIS}}(E_{\vert U})^{m-\textrm{cris}}_{n} \ar[r]^{1-\phi} \ar[d]^{{j''_{!}}(\beta)^m_{n}}  &  j_{!} f^{\ast}_{\textrm{CRIS}}(E_{\vert U})^{m-\textrm{cris}}_{n} \ar[r] \ar[d]^{{j''_{!}}(\beta)^m_{n}} & 0 \\
  (S_{2})\ 0 \ar[r]  & j''_{!}  \mathcal{G}_{n} \ar[r]& j''_{!} \mathcal{M}^{m-\textrm{cris}}_{n} \ar[r]^{1-\phi} &  j''_{!} \mathcal{M}^{m-\textrm{cris}}_{n} \ar[r]  & 0
}
$$

Le groupe $G$ agit de mani\`ere \'equivariante sur la suite exacte $(S_{1})$, donc sur le triangle distingu\'e $\mathcal{C}''(S_{1})$ obtenu en lui appliquant le foncteur

$$
\mathcal{C}'' := R \displaystyle \mathop{\lim}_{\leftarrow \atop{m}} \{(R  \displaystyle \mathop{\lim}_{\leftarrow \atop{n}} R \Gamma(\overline{X''}_{\textrm{SYNT},} -)) \otimes \mathbb{Q} \}.
$$\\

Comme $\mathcal{C}''(j''_{!}(\beta)^m_{n}) =: \tilde{\beta}$ est un isomorphisme, et de m\^eme en rempla\c{c}ant $\beta$ par $\alpha$, ou $\alpha_{\mathcal{G}}$, $\beta_{\mathcal{G}}$, le triangle distingu\'e $\mathcal{C}''(S_{2})$ obtenu en appliquant $\mathcal{C}''$ \`a $(S_{2})$ est aussi $G$-\'equivariant par transport de structure par ces isomorphismes : compte tenu de [E-LS 2, (3.11)] ce morphisme de triangles s'identifie \`a\\

$$
\xymatrix{
(3.3.13.1)\ R \Gamma_{\textrm{\'et},c}(U'',f^{\ast}(\mathcal{F}_{\vert U})) \otimes \mathbb{Q} \ar[r]  \ar[d] _{R \Gamma_{\textrm{\'et},c} (\beta_{\mathcal{G}}) \otimes \mathbb{Q}}^{\simeq} & R \Gamma_{\textrm{rig},c}(U'',f^{\ast}_{\textrm{rig}}(E^{\dag}_{K_{\vert U}})) \ar[r]^{1-\phi} \ar[d]_{R \Gamma_{\textrm{rig},c} (\beta^{\dag}_{K})}^{\simeq} & R \Gamma_{\textrm{rig},c}(U'',f^{\ast}_{\textrm{rig}}(E^{\dag}_{K_{\vert U}})) \ar[d]_{R \Gamma_{\textrm{rig},c} (\beta^{\dag}_{K})}^{\simeq}\\
(3.3.13.2)\qquad R \Gamma_{\textrm{\'et},c}(U'',\mathcal{G}) \otimes \mathbb{Q} \ar[r]  & R \Gamma_{\textrm{rig},c}(U'',\mathcal{M}^{\dag}_{K})  \ar[r]_{1-\phi} & R \Gamma_{\textrm{rig},c}(U'',\mathcal{M}^{\dag}_{K}) .
}
$$

En prenant les points fixes sous $G$ dans l'isomorphisme

$$
R \Gamma_{\textrm{rig},c}(U'', f^{\ast}_{\textrm{rig}}(E^{\dag}_{K \vert U})) \displaystyle \mathop{\longrightarrow}^{\sim} \mathcal{C}''(j''_{!} f^{\ast}_{\textrm{CRIS}}(E_{\vert U})^{m-\textrm{cris}}_{n})
$$

\noindent et en prouvant \`a la mani\`ere du th\'eor\`eme (3.2) que les points fixes sous $G$ du membre de droite s'identifient \`a
$$
\mathcal{C}(j_{U!}\ E_{\vert U_{n}}^{m-\textrm{cris}}) :=  R \displaystyle \mathop{\lim}_{\leftarrow \atop{m}} \{ (R  \displaystyle \mathop{\lim}_{\leftarrow \atop{n}} R \Gamma(\overline{X}_{\textrm{SYNT}}, j_{U!}\  E_{\vert U_{n}}^{m-\textrm{cris}})) \otimes \mathbb{Q} \} 
$$

\noindent on a prouv\'e le (1) du th\'eor\`eme (3.3.13) pour $U$ gr\^ace \`a [Et 3, IV, th\'eo (4.2)] ou [Et 6,  th\'eo (4.2)].\\

On a donc deux triangles distingu\'es reli\'es par des fl\`eches qui sont des isomorphismes

$$
\xymatrix{
R \Gamma_{\textrm{rig},c}(U, E^{\dag}_{K_{\vert U}}) \ar[r] \ar[d]^{\simeq} &  R \Gamma_{\textrm{rig},c}(X, E^{\dag}_{K }) \ar[r] &  R \Gamma_{\textrm{rig},c}(Z, E^{\dag}_{K_{\vert Z}}) \ar[d]^{\simeq} \\
\mathcal{C} (j_{U!}\ E_{\vert U_{n}}^{m-\textrm{cris}}) \ar[r] & \mathcal{C} (j_{!}\ E_{n}^{m-\textrm{cris}}) \ar[r] & \mathcal{C} (j_{Z!}\ E_{\vert Z_{n}}^{m-\textrm{cris}})\   ;
}
$$

\noindent d'apr\`es les axiomes des cat\'egories triangul\'ees [H, I, $\S1$] on peut compl\'eter par un isomorphisme au milieu. \\
Ceci ach\`eve la preuve du (1) du th\'eor\`eme.\\

\noindent \textit{Prouvons \`a pr\'esent le (2)}. Comme $k$ est alg\'ebriquement clos les points fixes sous $G$ de $R \Gamma_{\textrm{\'et},c}(U'', f^{\ast}(\mathcal{F}_{\vert U})) \otimes \mathbb{Q}$ sont \'egaux \`a $R \Gamma_{\textrm{\'et},c}(U, \mathcal{F}_{\vert U}) \otimes \mathbb{Q}$. Par suite les points fixes sous $G$ du triangle distingu\'e (3.3.13.1) $G$-\'equivariant fournissent une suite exacte longue de cohomologie\\

\noindent (3.3.13.3) $ \rightarrow H^i_{\textrm{\'et},c}(U, \mathcal{F}_{\vert U}) \otimes \mathbb{Q} \rightarrow H^i_{\textrm{rig},c}(U/K, E^{\dag}_{K_{\vert U}}) \displaystyle \mathop{\rightarrow}_{1-\phi} H^i_{\textrm{rig},c}(U/K, E^{\dag}_{K_{\vert U}}) \rightarrow .$\\

\noindent Les groupes de cohomologie rigide $H^i_{\textrm{rig},c}(U/K, E^{\dag}_{K_{\vert U}})$ \'etant de dimension finie sur $K$ [Tsu 1, theo 6.1.2], et $k$ alg\'ebriquement clos, ces suites exactes se scindent en suites exactes courtes [I$\ell$; II, lemme 5.6]. Notons $Z \hookrightarrow X$ l'immersion du ferm\'e compl\'ementaire \`a $U$ (rappelons que $X$ est suppos\'e r\'eduit). Si dim $Z = 0$, $Z$ est fini \'etale sur $k$ et on a une suite analogue \`a (3.3.13.3) pour $Z$ : par r\'ecurrence sur la dimension de $X$ on peut donc supposer l'existence de suites exactes courtes telles que (3.3.13.3) pour $U$ et pour $Z$. En particulier on a trois triangles distingu\'es horizontaux reli\'es par des fl\`eches induites par (3.3.13.3) appliqu\'e \`a $U$ et $Z$ :

$$
\begin{array}{c}
\xymatrix{
R \Gamma_{\textrm{\'et},c}(U, \mathcal{F}_{U \mathbb{Q}}) \ar[r] \ar[d] & R \Gamma_{\textrm{\'et},c}(X, \mathcal{F}_{ \mathbb{Q}}) \ar[r]  & R \Gamma_{\textrm{\'et},c}(Z, \mathcal{F}_{Z \mathbb{Q}}) \ar[d]\\
R \Gamma_{\textrm{rig},c}(U,  E^{\dag}_{K_{\vert U}}) \ar[r] \ar[d]_{1-\phi_{U}} & R \Gamma_{\textrm{rig},c}(X, E^{\dag}_{K}) \ar[r] \ar[d]_{1-\phi} & R \Gamma_{\textrm{rig},c}(Z, E_{K_{\vert Z}}^{\dag}) \ar[d]^{1-\phi_{Z}} \\
R \Gamma_{\textrm{rig},c}(U,  E^{\dag}_{K_{\vert U}}) \ar[r]  & R \Gamma_{\textrm{rig},c}(X, E^{\dag}_{K}) \ar[r]  & R \Gamma_{\textrm{rig},c}(Z, E_{K_{\vert Z}}^{\dag}) ; 
}
\end{array}
\leqno{(3.3.13.4)}
$$

\noindent par les axiomes des cat\'egories triangul\'ees [H, I, $\S$\ 1] on peut compl\'eter par un morphisme $R \Gamma_{\textrm{\'et},c} (X, \mathcal{F}_{\mathbb{Q}}) \rightarrow R \Gamma_{\textrm{rig},c}(X, E^{\dag}_{K})$.
\noindent Les suites exactes courtes (3.3.13.3) pour $U$ et $Z$ fournissent alors l'analogue pour $X$, d'o\`u le (2) du th\'eor\`eme (3.3.13).\\

\noindent \textit{Autre d\'emonstration du (2)}. Une autre m\'ethode consiste \`a appliquer le foncteur $\mathcal{C}$ \`a la suite exacte du th\'eor\`eme (3.3.12)\\

\noindent (3.3.13.5) $ 0 \rightarrow j_{!}\  \mathcal{F}_{n} \rightarrow j_{!}\ (\mathcal{F}_{n} \otimes_{\mathcal{V}^{\sigma}_{n}} \mathcal{O}^{m-\textrm{cris}}_{n,X}) \displaystyle \mathop{\rightarrow}_{1-\phi}  j_{!}\ (\mathcal{F}_{n} \otimes_{\mathcal{V}^{\sigma}_{n}} \mathcal{O}^{m-\textrm{cris}}_{n,X}) \rightarrow 0$

\noindent o\`u l'on remarque que $E^{m-\textrm{cris}}_{n} =  \mathcal{F}_{n} \otimes_{\mathcal{V}^{\sigma}_{n}} \mathcal{O}^{m-\textrm{cris}}_{n,X}$. Par le (1) du th\'eor\`eme le triangle distingu\'e ainsi obtenu s'identifie au triangle distingu\'e

$$
R \Gamma_{\textrm{\'et},c} (X, \mathcal{F}_{\mathbb{Q}}) \longrightarrow R \Gamma_{\textrm{rig},c}(X, E^{\dag}_{K}) \displaystyle \mathop{\rightarrow}_{1-\phi} R \Gamma_{\textrm{rig},c}(X, E^{\dag}_{K}).
$$

\noindent La suite exacte longue de cohomologie se scinde alors en suites exactes courtes par le m\^eme argument que ci-dessus. $\square$

\vskip 3mm
\noindent \textbf{Remarque (3.3.14)}. En supposant seulement que $k$ contient $\mathbb{F}_{q}$ [cf (3.0)] et que $\mathcal{F}_{n}$ est un $\mathcal{V}^{\sigma}_{n}$-module localement trivial sur SYNT$(X)$, on pose encore

$$
E^{m-\textrm{cris}}_{n} = T^{(m)}(\mathcal{F}_{n})^{\textrm{cris}} = \mathcal{F}_{n} \otimes_{\mathcal{V}^{\sigma}_{n}} \mathcal{O}^{m-\textrm{cris}}_{n,X}.
$$

\noindent En appliquant le foncteur $R \Gamma(\overline{X}_{\textrm{SYNT}, -})$ \`a la suite exacte (3.3.13.5), on obtient un triangle distingu\'e

$$
R \Gamma_{\textrm{\'et},c} (X, \mathcal{F}_{n}) \longrightarrow R \Gamma_{\textrm{synt},c}(X, E^{m-\textrm{cris}}_{n}) \displaystyle \mathop{\longrightarrow}_{1-\phi} R \Gamma_{\textrm{synt},c}(X, E^{m-\textrm{cris}}_{n}).
$$

\end{document}